\documentclass{amsproc} 
\usepackage[mathscr]{eucal}
\usepackage{amssymb}
\usepackage{latexsym}
\usepackage{amsmath,amsthm,amscd}
\font\es=eufm10

\def\gge{\mbox{\es {e}}} 
 
\def\gg{\mbox{\es {g}}}

\def\gC{\mbox{\es {C}}}

\def\gJ{\mbox{\es {J}}}

\def\gP{\mbox{\es {P}}}

\def\diag{\mbox{\rm {diag}}}

\def\Ad{\mbox{\rm {Ad}}}

\def\Iso{\mbox{\rm {Iso}}}

\def\ov{\overline}
\def\wti{\widetilde}

\def\C{\mbox{\boldmath $C$}}

\def\R{\mbox{\boldmath $R$}}

\def\Z{\mbox{\boldmath $Z$}}

\def\sR{\mbox{\boldmath $\scriptscriptstyle{R}$}}
\def\sC{\mbox{\boldmath $\scriptscriptstyle{C}$}}

\def\0{\mbox{\boldmath {0}}}    
\def\1{\mbox{\boldmath {1}}}      
\def\2{\mbox{\boldmath {2}}}      
\def\3{\mbox{\boldmath {3}}}      
\def\4{\mbox{\boldmath {4}}}      
\def\5{\mbox{\boldmath {5}}}      
\def\6{\mbox{\boldmath {6}}}      
\def\7{\mbox{\boldmath {7}}}      
\def\8{\mbox{\boldmath {8}}}      
\def\9{\mbox{\boldmath {9}}}

\def\m{\mbox{\boldmath $m$}}

\begin{document}

\title[Decomposition of compact exceptional Lie groups ]
{Decomposition of compact exceptional Lie groups \\into their maximal tori}

\author[Toshikazu MIYASHITA]{Toshikazu MIYASHITA}

\maketitle


\thispagestyle{empty}
\vspace{-6mm}

\begin{center}
{\footnotesize \it Nagano prefectural Komoro high school}
\end{center}
\vspace{4mm}

{\small {\bf Abstract.}\quad In this paper we treat the intersection of fixed point subgroups by 
the involutive automorphisms of exceptional Lie group $G= F_4,\ E_6, E_7$. 
We shall find involutive automorphisms of $G$ such that the connected 
component of the intersection
 of those fixed point subgroups coincides with the maximal torus of $G$. }
\vspace{6mm}

{\bf 1. Introduction}
\vspace{2mm}

 It is known that the involutive automorphisms of the compact 
 Lie groups play an important role in the theory of symmetric space 
 (c.f. Berger [1]). In [8],[9] Yokota showed that the exceptional symmetric 
 spaces $G/H$ are realized definitely by calculating the fixed point subgroup of 
the involutive automorphisms ${\tilde \gamma},  {\tilde \gamma'}, {\tilde \sigma},{\tilde \sigma'},{\tilde \iota}$
 of $G$, where ${\tilde \gamma}, {\tilde \gamma'}, {\tilde \sigma},{\tilde \sigma'}$ are induced by $\R$-linear transformations $\gamma, \gamma', \sigma, \sigma'$ of $\gJ$ and ${\tilde \iota}$ is induced by $C$-linear transformation $\iota$ of $\gP^C$. Here $\gamma, \gamma' \in G_2 \subset F_4 \subset E_6 \subset 
 E_7, \sigma, \sigma' \in F_4 \subset E_6 \subset E_7$ and $\iota \in E_7$. For the cases of the graded 
 Lie algebras $\gg$ of the second kind and third kind, the corresponding subalgebras
  ${\gg}_0,\ {\gg}_{ev},\ {\gg}_{ed}$ 
of $\gg$ are realized as the intersection of those fixed point subgroups of 
the commutative involutive automorphisms ([3],[6],[7],[10],[11],[12]). \par 

  In [2],[4],[5] we determined the intersection of those fixed point subgroups of 
the involutive automorphisms of $G$ when $G$ is a compact exceptional Lie group. 
We remark that those intersection subgroups are maximal rank of $G$. \par

  In general, let $G$ be a connected compact Lie group and $\sigma_1, \sigma_2, \cdots, 
\sigma_m$ commutative automorphism elements of $G$. Set $G^{\sigma_1, 
\sigma_2, \cdots, \sigma_k} = \{\alpha \in G \, | \, \sigma_i\alpha = 
\alpha\sigma_i, i = 1, \cdots, k\}$. We expect that the group $G^{\sigma_1, 
\sigma_2, \cdots, \sigma_k}$ is a maximal rank subgroup of $G$. Consider the following degreasing sequence of 
subgroups of $G$:
$$
   G^{\sigma_1} \supset G^{\sigma_1, \sigma_2} \supset \cdots \supset 
G^{\sigma_1, \cdots, \sigma_m}. $$
Let $T^l$ be the maximal tours of $G$.
In this paper we would like to find $\sigma_1, \sigma_2, \cdots, \sigma_m$ such that 
the connected component subgroup 
$(G^{\sigma_1, \sigma_2, \cdots \sigma_k})_0$ of the group $G^{\sigma_1, \sigma_2, \cdots \sigma_k}$ is isomorphic to $T^l$ when $G$ is simply connected compact exceptional Lie groups $G_2, F_4, E_6$ or $E_7$. For the case $G=G_2$, we prove that the group $((G_2)^{\gamma, \gamma'})_0 \cong T^2$ by [5], Theorem 1.1.3.
Then we shall prove the following\,:

$$
\begin{array}{ll}
\hspace*{-10mm} (1) &  ((F_4)^{\gamma, \gamma', \sigma, \sigma'})_0 \cong T^4 , 
\vspace{2mm}\\
\hspace*{-10mm} (2) &  ((E_6)^{\gamma, \gamma', \sigma, \sigma'})_0 \cong T^6 , 
\vspace{2mm}\\
\hspace*{-10mm} (3) &  ((E_7)^{\gamma, \gamma', \sigma, \sigma', \iota})_0 \cong T^7 .
\end{array} 
$$
For the case $G=E_8$, we conjecture that the group $((E_8)^{\gamma, \gamma', \sigma, 
\sigma',\upsilon_3})_0 \cong T^8$, where $\lambda' \in E_8$ (As for $\upsilon_3$, see [3]).
\vspace{4mm}

{\bf 2. Group $F_4$}
\vspace{2mm}

The simply connected compact Lie group $F_4$ is given by the automorphism group of the exceptional Freudenthal algebra $\gJ$\,:
$$
   F_4 = \{\alpha \in \Iso_{\sR}(\gJ) \, | \, \alpha(X \times Y) = \alpha X 
\times \alpha Y \}. $$

We shall review the definitions of $\R$-linear transformations $\gamma, \gamma', \sigma, \sigma'$ of $\gJ$([8], [10], [12]).
\vspace{2mm}

Firstly we define $\R$-linear transformations $\gamma, \gamma'$ and $\gamma_1$ of $\gJ_{\sC} \oplus M(3, \C) = \gJ$ by
\begin{eqnarray*}     
    \gamma(X + M) \!\!\!&=&\!\!\! X + \gamma(\m_1, \m_2, \m_3) = X + (\gamma\m_1, \gamma\m_2, \gamma\m_3),
\vspace{1mm}\\
   \gamma'(X + M) \!\!\!&=&\!\!\! X + \gamma'(\m_1, \m_2, \m_3) = X + (\gamma'\m_1, \gamma'\m_2, \gamma'\m_3), 
\vspace{1mm}\\
     \gamma_1(X + M) \!\!\!&=&\!\!\! \ov{X} + \ov{M},\quad X + M \in \gJ_{\sC} \oplus M(3, \C) = \gJ,  
\end{eqnarray*}
respectively, where $\gJ_{\sC}=\{X \in M(3, \C)\,|\,X^* =X \}$, \vspace{1mm}the right-hand side transformations $\gamma, \gamma' : \C^3 \to 
\C^3$ are defined by 
$$
\begin{array}{c}
     \gamma\Big(\begin{pmatrix}n_1 \cr
                         n_2 \cr
                         n_3
                      \end{pmatrix}\Big) =  \begin{pmatrix}n_1 \cr
                         -n_2 \cr
                         -n_3
                        \end{pmatrix}, \quad
     \gamma'\Big(\begin{pmatrix}n_1 \cr
                         n_2 \cr
                         n_3
                          \end{pmatrix}\Big) = \begin{pmatrix}-n_1 \cr
                         n_2 \cr
                         -n_3
                         \end{pmatrix},\,\,n_i \in \C.
\end{array} 
$$ 
Then $\gamma, \gamma', \gamma_1 \in G_2 \subset F_4$, and $\gamma^2={\gamma'}^2={\gamma_1}^2=1$.
\vspace{2mm}

Further we define $\R$-linear transfomations $\sigma$ and 
$\sigma'$ of $\gJ_{\sC} \oplus M(3, \C) = \gJ$ by
\begin{eqnarray*}
   \sigma(X + M)\!\!\!&=&\!\!\! \sigma X + (\m_1, -\m_2, -\m_3),
\vspace{1mm}\\
   \sigma'(X + M)\!\!\!&=&\!\!\! \sigma' X + (-\m_1, -\m_2, \m_3), \quad X + M \in \gJ_{\sC} \oplus M(3, \C) = \gJ,  
\end{eqnarray*}
respctively, where the right-hand side transformations $\sigma, \sigma' : \gJ_{\sC} \to 
\gJ_{\sC}$ are defined by 
$$
\sigma X=\sigma \begin{pmatrix} \xi_1     & x_3      & \ov{x}_2 \\
                                \ov{x}_3  & \xi_2    & x_1      \\
                                x_2       & \ov{x}_1 & \xi_3
                \end{pmatrix}
               =\begin{pmatrix} \xi_1     & -x_3     & -\ov{x}_2 \\
                                -\ov{x}_3 & \xi_2    & x_1      \\
                                -x_2      & \ov{x}_1 & \xi_3
                \end{pmatrix},\,\,
\sigma' X=     \begin{pmatrix} \xi_1     & x_3       & -\ov{x}_2 \\
                                \ov{x}_3 & \xi_2     & -x_1      \\
                                -x_2     & -\ov{x}_1 & \xi_3
                \end{pmatrix}.               
$$
Then $\sigma, {\sigma'} \in F_4$ and $\sigma^2 = {\sigma'}^2 =  
1$.
\vspace{2mm}

The group $\Z_2 = \{1, \gamma_1\}$ acts on the group $U(1) \times U(1) \times 
SU(3)$ by
$$
         \gamma_1(p, q, A) = (\ov{p}, \ov{q}, \ov{A}). 
$$
Hence the group $\Z_2 = \{1, \gamma_1\}$ acts naturally on the group $(U(1) \times U(1) \times SU(3))/\Z_3$. 
\vspace{1mm}

Let $(U(1) \times U(1) \times SU(3)) \cdot \Z_2$ be the semi-direct product of 
those groups under this action. 
\vspace{1mm}

Hereafter, $\omega_1$ denotes $-\dfrac{1}{2}+\dfrac{\sqrt{3}}{2} e_1 \in \gC.$
\vspace{2mm}

{\sc Proposition 2.1.} $(F_4)^{\gamma,\gamma'} \cong ((U(1) \times U(1) \times 
SU(3))/\Z_3) \cdot \Z_2, \; \vspace{0.5mm} \Z_3 = \{(1, 1, E),$ $ (\omega_1, \omega_1, 
{\omega_1}E), ({\omega_1}^2, {\omega_1}^2, {\omega_1}^2 E)\}$.
\vspace{2mm}

{\sc Proof.} We define a mapping $\varphi_4 : (U(1) \times U(1) \times SU(3)) \cdot \Z_2 \to (F_4)^{\gamma, \gamma'}$ by
$$
\begin{array}{l}
     \;\,\varphi_4((p, q, A), 1)(X + M) = AXA^* + D(p, q)MA^*, 
\vspace{1.5mm}\\
     \varphi_4((p, q, A), \gamma_1)(X + M) = A\ov{X}A^* + D(p, q)\ov{M}A^*, 
\vspace{1.5mm}\\
\hspace{6cm}
     X + M \in \gJ_{\sC} \oplus M(3, \C) = \gJ,
\vspace{-2mm}
\end{array} $$
where $D(p, q) = \diag(p, q, \ov{pq}) \in SU(3)$. Then $\varphi_4$ induces the 
required isomorphism (see [5] for details). \;\; \hfill $\Box$
\vspace{2mm}

{\sc Lemma 2.2.}\,\,{\it  The mapping $\varphi_4 : (U(1) \times U(1) \times SU(3)) \cdot \Z_2 \to {(F_4)}^{\gamma, \gamma'}$ satisfies 
$$
    \sigma  = \varphi_4((1, 1, E_{1,-1}),1),\,\,\,\, 
    \sigma' = \varphi_4((1, 1, E_{-1,1}),1), $$
where $E_{1,-1}= \diag(1, -1, -1),\,\, E_{-1, 1} = \diag(-1, -1, 1) \in SU(3)$}.
\vspace{2mm}

We denote $U(1) \times  \cdots \times U(1), (1, \cdots 1)$ and $(\omega_k, \cdots \omega_k)$ ($l$-times) by $U(1)^{\times l}, (1)^{\times l}$ and $(\omega_k)^{\times l}$,  
respectively.  
\vspace{2mm}

Now, we determine the structures of the group ${(F_4)}^{\gamma, \gamma', 
\sigma, \sigma'}= ((F_4)^{\gamma, \gamma'})^{\sigma, \sigma'}$. 
\vspace{3mm}

{\sc Theorem 2.3.} \quad \qquad \qquad $((F_4)^{\gamma,\gamma',\sigma,\sigma'})_0 \cong U(1)^{\times 4}$. 
\vspace{2mm}

{\sc Proof.} For $\alpha \in (F_4)^{\gamma, \gamma', \sigma, \sigma'} \subset
(F_4)^{\gamma, \gamma'}$, there exist $p, q \in U(1)$ and $A \in SU(3)$ such 
that $\alpha = \varphi_4((p, q, A), 1)$ or $\alpha = \varphi_4((p, q, A), 
\gamma_1)$ (Proposition 2.1). For the case of $\alpha=\varphi_4((p, q, A), 1)$, by combining  the conditions of 
$\sigma\alpha\sigma = \alpha$ and $\sigma'\alpha\sigma' = \alpha$ with Lemma 2.2, we have 
$$
    \varphi_4((p, q, E_{1,-1} A E_{1,-1}),1) = \varphi_4((p, q, A), 1) $$
\vspace{-1mm} 
and
\vspace{-1mm}
$$
  \varphi_4((p, q, E_{-1,1} A E_{-1,1}),1) = \varphi_4((p, q, A), 1) . $$
Hence
$$
\text{(i)}\,\,       E_{1,-1} A E_{1, -1} = A,
\,\,\,\,\text{(ii)}\left\{ \begin{array}{l}
       p = {\omega_1}p  
\vspace{1mm}\\
       q = {\omega_1}q  
\vspace{1mm}\\
       E_{1,-1} A E_{1, -1} = {\omega_1}A,
       \end{array} \right.
\,\,\,\,\text{(iii)}\left\{ \begin{array}{l}
       p = {\omega_1}^2 p  
\vspace{1mm}\\
       q = {\omega_1}^2 q  
\vspace{1mm}\\
       E_{1,-1} A E_{1, -1} = {\omega_1}^2 A
       \end{array} \right.
$$
and
$$
\text{(iv)}\,\,       E_{-1,1} A E_{-1, 1} = A,
\,\,\,\,\text{(v)}\left\{ \begin{array}{l}
       p = {\omega_1}p  
\vspace{1mm}\\
       q = {\omega_1}q  
\vspace{1mm}\\
       E_{-1,1} A E_{-1, 1} = {\omega_1}A,
       \end{array} \right.
\,\,\,\,\text{(vi)}\left\{ \begin{array}{l}
       p = {\omega_1}^2p  
\vspace{1mm}\\
       q = {\omega_1}^2q  
\vspace{1mm}\\
       E_{-1,1} A E_{-1, 1} = {\omega_1}^2A.
       \end{array} \right.
$$
We can eliminate the case (ii), (iii), (v) or (vi) because $p \ne 0$ or $q \ne 0$.
Hence we have $p, q \in U(1)$ and $ A \in S(U(1) \times U(1) \times 
U(1))$. Since the mapping 
$U(1) \times U(1)  \to S(U(1) \times U(1) \times U(1))$, 
$$
     h(a_1, a_2) = (a_1, a_2, \ov{a_1a_2}) $$
is an isomorphism, the group satisfying with the conditions of case (i) and (iv) is $(U(1)^{\times 4})/\Z_3$. 
For the case of $\alpha = \varphi_4((p, q, A), \gamma_1)$, from $\varphi_4((p, q, A),
\gamma_1)\! =\! \varphi_4((p, q, A), $ $1)\gamma_1, \varphi_4((1, 1, 
E_{1,-1}), 1)\gamma_1=\gamma_1 \varphi_4((1, 1, E_{1, -1}), 1)$ and 
$\varphi_4((1, 1, E_{-1,1}), 1)\gamma_1\!=\!\gamma_1 \varphi_4((1,$ $ 1, E_{-1, 1}), 
1)$ , this case is in the same situation as above. \vspace{0.5mm}Thus we have 
$(F_4)^{\gamma,\gamma',\sigma,\sigma'} \!\cong\! ((U(1)^{\times 4})/\Z_3) \cdot 
\Z_2, \Z_3 \!=\! \{(1)^{\times 4}, (w_1)^{\times 4}, ({w_1}^2)^{\times 4}\}$. 
\vspace{0.5mm}The group $(U(1)^{\times 4})/\Z_3$ is naturally isomorphic to the torus 
$U(1)^{\times 4}$, hence we obtain $(F_4)^{\gamma,\gamma',\sigma,\sigma'} \cong (U(1)^{\times 4}) \cdot \Z_2$. Therefore we have the required isomorphism of the theorem. \hfill $\Box$
\vspace{3mm}

{\bf 3. The group $E_6$}
\vspace{2mm}

The simply connected compact Lie group $E_6$ is given by
$$
    E_6 = \{\alpha \in \Iso_C(\gJ^C) \, | \, \alpha X \times \alpha Y = 
\tau\alpha\tau(X \times Y), \langle \alpha X, \alpha Y \rangle = \langle X, Y 
\rangle\}. $$

$\R$-linear transformations $\gamma, \gamma', \gamma_1, \sigma$ and $\sigma'$ 
of $\gJ = \gJ_{\sC} \oplus M(3, \C)$ \vspace{0.5mm}are naturally extended to the $C$-linear 
transformations of $\gamma, \gamma', \gamma_1, \sigma$ and $\sigma'$ of $\gJ^C 
= (\gJ_{\sC})^C \oplus M(3, \C)^C$. Then we have $\gamma, \gamma', \gamma_1, \sigma, \sigma' \in E_6$.
\vspace{2mm}

The group $\Z_2 = \{1, \gamma_1\}$ acts on the group $U(1) \times U(1) \times 
SU(3) \times SU(3)$ by
$$
    \gamma_1(p, q, A, B) = (\ov{p}, \ov{q}, \ov{B}, \ov{A}). 
$$
Hence the group $\Z_2 = \{1, \gamma_1\}$ acts naturally on the group $(U(1) \times U(1) \times SU(3) \times SU(3))/\Z_3$.
\vspace{2mm}

Let $(U(1) \times U(1) \times SU(3) \times SU(3)) \cdot \Z_2$ be the semi-direct product of those groups under this action.
\vspace{3mm}

{\sc Proposition 3.1.} $(E_6)^{\gamma, \gamma'} \cong ((U(1) \times U(1) \times SU(3) \times SU(3))/\Z_3) \cdot \Z_2, \vspace{0.5mm}\Z_3 = \{(1, 1, E, E), (\omega_1, 
\omega_1, \omega_1E, \omega_1E), ({\omega_1}^2, {\omega_1}^2, {\omega_1}^2E, 
{\omega_1}^2E)\}$.
\vspace{2mm}

{\sc Proof.} We define a mapping $\varphi_6 : (U(1) \times U(1) \times SU(3) 
\times SU(3)) \cdot \Z_2 \to (E_6)^{\gamma, \gamma'}$ by
$$
\begin{array}{lll}
    \varphi_6((p, q, A, B), 1)(X + M) = h(A, B)Xh(A, B)^* + D(p, q)M\tau h(A, B)^*,
\vspace{1.5mm}\\
    \varphi_6((p, q, A, B), \gamma_1)(X + M) = h(A, B)\ov{X}h(A, B)^* + D(p, q)\ov{M}\tau h(A, B)^*, 
\vspace{1.5mm}\\
\hspace{6cm} X + M \in (\gJ_{\sC})^C \oplus M(3, \C)^C = \gJ^C.
\end{array} $$
Here $D(p, q) = \diag(p, q, \ov{pq}) \in SU(3)$ and $h : M(3, \C) \times M(3, 
\C) \to M(6, \C)^C$ is defined by
$$
       h(A, B) = \dfrac{A + B}{2} + i\dfrac{A - B}{2}e_1. $$
Then $\varphi_6$ induces the required isomorphism (see [5] for details).\hfill $\Box$
\vspace{3mm}

{\sc Lemma 3.2.} {\it The mapping $\varphi_6 : (U(1) \times U(1) \times SU(3) \times SU(3)) \cdot \Z_2 \to (E_6)^{\gamma, \gamma'}$ satisfies}
$$
    \sigma = \varphi_6((1, 1, E_{1, -1}, E_{1,-1}), 1), \quad 
    \sigma' = \varphi_6((1, 1, E_{-1, 1}, E_{-1, 1}), 1). $$ 
\vspace{-1mm}

The group $\Z_2 = \{1, \gamma_1\}$ acts on the group $U(1)^{\times 6}$ by
$$
    \gamma_1(p, q, a_1, a_2, a_3, a_4) = (\ov{p}, \ov{q}, \ov{a}_3, \ov{a}_4, \ov{a}_1, \ov{a}_2). $$
Let $(U(1)^{\times 6}) \cdot \Z_2$ be the semi-direct product of those groups under this action.
\vspace{2mm}

Now, we determine the structures of the gruop $(E_6)^{\gamma, \gamma', \sigma, \sigma'}= ((E_6)^{\gamma, \gamma'})^{\sigma, \sigma'}$. 
\vspace{2mm}

{\sc Theorem 3.3.} \qquad \qquad $((E_6)^{\gamma, \gamma', \sigma, \sigma'})_0 \cong U(1)^{\times 6}$.  
\vspace{2mm}

{\sc Proof.} For $\alpha \in (E_6)^{\gamma, \gamma', \sigma, \sigma'} \subset 
(E_6)^{\gamma, \gamma'}$, there exist $p, q \in U(1)$ and $A, B \in SU(6)$ such that $\alpha = \varphi_6((p, q, A, B), 1)$ or $\alpha = \varphi_6((p, q, A, B), \gamma_1)$ (Proposition 3.1). For the case of $\alpha = \varphi_6((p, q, A, B), 1)$, by combining the conditions $\sigma\alpha\sigma = \alpha$ and $\sigma'\alpha\sigma' = \alpha$ with Lemma 3.2, we have
$$
    \varphi_6((p, q, E_{1, -1}AE_{1, -1}, E_{1, -1}BE_{1, -1}), 1) = \varphi_6((p, q, A, B), 1) $$
\vspace{-1mm}
and
\vspace{-1mm}
$$
    \varphi_6((p, q, E_{-1, 1}AE_{-1, 1}, E_{-1, 1}BE_{-1, 1}), 1) = \varphi_6((p, q, A, B), 1). $$
Hence
$$
\text{(i)}\left\{\begin{array}{l}
  E_{1, -1}AE_{1, -1} = A
\vspace{1mm}\\
   E_{1, -1}BE_{1, -1} = B,
\end{array} \right.  
\text{(ii)}\left\{\begin{array}{l}
   p = \omega_1p
\vspace{1mm}\\
   q = \omega_1q
\vspace{1mm}\\
   E_{1, -1}AE_{1, -1} = \omega_1A
\vspace{1mm}\\
   E_{1, -1}BE_{1, -1} = \omega_1B,
\end{array} \right.
\text{(iii)}\left\{\begin{array}{l}
   p = {\omega_1}^2p
\vspace{1mm}\\
   q = {\omega_1}^2q
\vspace{1mm}\\
   E_{1, -1}AE_{1, -1} = {\omega_1}^2A
\vspace{1mm}\\
   E_{1, -1}BE_{1, -1} = {\omega_1}^2B
\end{array} \right. $$
and
$$
\text{(iv)}\left\{\begin{array}{l}
   E_{-1, 1}AE_{-1, 1} = A
\vspace{1mm}\\
   E_{-1, 1}BE_{-1, 1} = B,
\end{array} \right.  
\!\text{(v)}\left\{\begin{array}{l}
   p = \omega_1p
\vspace{1mm}\\
   q = \omega_1q
\vspace{1mm}\\
   E_{-1, 1}AE_{-1, 1} = \omega_1A
\vspace{1mm}\\
   E_{-1, 1}BE_{-1, 1} = \omega_1B,
\end{array} \right.
\!\text{(vi)}\left\{\begin{array}{l}
   p = {\omega_1}^2p
\vspace{1mm}\\
   q = {\omega_1}^2q
\vspace{1mm}\\
   E_{-1, 1}AE_{-1, 1} = {\omega_1}^2A
\vspace{1mm}\\
   E_{-1, 1}BE_{-1, 1} = {\omega_1}^2B.
\end{array} \right. $$
We can eliminate the case (ii), (iii), (v) or (vi) because $p \ne 0$ or $q \ne 0$. Thus we have $p, q \in U(1)$ and $A, B \in S(U(1)^{\times 3})$. Since the 
mapping 
$U(1)^{\times 4} \to S(U(1)^{\times 5})$, 
$$
     h(a_1, a_2, a_3, a_4) = (a_1, a_2, a_3, a_4,\ov{a_1a_2a_3a_4}) $$
is an isomorphism, the group satisfying with the conditions of case (i) and (iv) is $(U(1)^{\times 6})/\Z_3$. For the case of $\alpha=\varphi_6((p, q, A, B), \gamma_1)$, from $\varphi_6((p, q, A, B), \gamma_1) = 
\varphi_6((p, q,A,B),1)\gamma_1,\, \varphi_6((1, 1, E_{1, -1}, E_{1, -1}), 1)\gamma_1\,=\, \gamma_1\varphi_6((1, 1,E_{1, -1}, E_{1, -1}), 1)$ and $\varphi_6((1,1, E_{-1, 1}, E_{-1, 1}), 1)\gamma_1 
= \gamma_1\varphi_6((1, 1, E_{-1, 1}, E_{-1, 1}), 1)$, this case is in the same situation as above. Thus we have 
$(E_6)^{\gamma,\gamma',\sigma,\sigma'} \!\!\cong ((U(1)^{\times 6})/\Z_3) \cdot 
\Z_2, \Z_3 \!=\! \{(1)^{\times 6}, (w_1)^{\times 6}, $ $({{w_1}^2})^{\times 6}\}$. The 
group $(U(1)^{\times 6})/\Z_3$ is naturally isomorphic to the torus 
$U(1)^{\times 6}$, hence we obtain $(E_6)^{\gamma, \gamma', \sigma, \sigma'} \cong (U(1)^{\times 6}) \cdot \Z_2$. Therefore we have the required isomorphism of the theorem.\hfill $\Box$
\vspace{3mm}

{\bf 4. Group $E_7$}
\vspace{2mm}

Let $\gP^C = \gJ^C \oplus \gJ^C \oplus C \oplus C$. The simply connected compact Lie group $E_7$ is given by
$$
    E_7 = \{\alpha \in \Iso_C(\gP^C) \, | \, \alpha(P \times Q)\alpha^{-1} = 
\alpha P \times \alpha Q, \langle \alpha P, \alpha Q \rangle = \langle P, Q 
\rangle\}. 
$$

Under the identification $(\gP_{\sC})^{C} \oplus (M(3, \C)^{C} \oplus M(3, 
\C)^{C})$ \vspace{0.5mm}with $\gP^{C}\;$: $((X, Y, \xi, \eta),$ $ (M, N)) = (X + M, Y + 
N, \xi, \eta)$, $C$-linear transformations of $\gamma, \gamma',\gamma_1, 
\sigma$ and $\sigma'$ of $\gJ^C$ are extended to $C$-linear transformations of 
$\gP^{C}$ as
\begin{eqnarray*}
    \; \gamma(X + M, Y + N, \xi, \eta) \!\!\!&=&\!\!\! (X + \gamma M, Y + 
\gamma N, \xi, \eta), 
\vspace{1mm}\\
    \, \gamma'(X + M, Y + N, \xi, \eta) \!\!\!&=&\!\!\! (X + \gamma'M, Y + 
\gamma'N, \xi, \eta),
\vspace{1mm}\\
    \gamma_1(X + M, Y + N, \xi, \eta) \!\!\!&=&\!\!\! (\ov{X} + \ov{M}, \ov{Y} 
+ \ov{N}, \xi, \eta),
\vspace{1mm}\\
\; \sigma(X + M, Y + N, \xi, \eta) \!\!\!&=&\!\!\!(\sigma X + \sigma M, 
\sigma 
Y + \sigma N, \xi, \eta), 
\vspace{1mm}\\
   \gamma(X + M, Y + N, \xi, \eta) \!\!\!&=&\!\!\!(\sigma' X + \sigma' M, 
\sigma' Y + \sigma' N, \xi, \eta),
\end{eqnarray*} 
where $\gamma M = \diag(1, -1, -1) M, \gamma' M = \diag(-1, -1, 1) M, \sigma M 
= M \diag(1, -1, -1)$ and $ \sigma'M = M \diag(-1, -1, 1)$.

\noindent Moreover we define a $C$-linear transformation $\iota$ of $\gP^C$ by
$$
  \iota((X + M, Y + N, \xi, \eta)=(-iX-iM, iY+iN, -i\xi, i\eta).
$$
\vspace{-3mm}

The group $\Z_2 = \{1, \gamma_1\}$ acts the group $U(1) \times U(1) \times 
SU(6)$ by
$$
    \gamma_1(p, q, A) = (\ov{p}, \ov{q}, \ov{(\Ad J_3)A)}, \quad J_3 = 
\begin{pmatrix} 0 & E \cr
                - E & 0 \end{pmatrix}. 
$$
Hence the group $\Z_2 = \{1, \gamma_1\}$ acts naturally on the group $(U(1) \times U(1) \times SU(6))/\Z_3$.
\vspace{-1mm}

Let $(U(1) \times U(1) \times SU(6)) \cdot \Z_2$ be the semi-direct product of those groups under this action.
\vspace{3mm}

{\sc Proposition 4.1.} $(E_7)^{\gamma, \gamma'} \cong ((U(1) \times U(1) \times SU(6))/\Z_3) \cdot \Z_2,\,\vspace{0.5mm} \Z_3 = \{(1, 1, E),$ $ (\omega_1,\omega_1, \omega_1E), ({\omega_1}^2, {\omega_1}^2, {\omega_1}^2E)\}$.
\vspace{2mm}

{\sc Proof.} We define a mapping $\varphi_7 : (U(1) \times U(1) \times SU(6)) \cdot \Z_2 \to (E_7)^{\gamma, \gamma'}$ by
\begin{eqnarray*}
    \varphi_7((p, q, A), 1)P \!\!\!&=&\!\!\! f^{-1}((D(p, q), A)(fP)),
\vspace{1.5mm}\\
    \varphi_7((p, q, A), \gamma_1)P \!\!\!&=&\!\!\! f^{-1}((D(p, q), A)(f\gamma_1P)), \quad P \in \gP^C.
\end{eqnarray*}
Here $D(p, q) = \diag(p, q, \ov{pq}) \in SU(3)$ and the mapping $f$ is defined in [9], Section 2.4. Then $\varphi_7$ induces the required isomorphism (see [5] for details). \hfill $\Box$
\vspace{3mm}

{\sc Lemma 4.2.} {\it The mapping $\varphi_7 : (U(1) \times U(1) \times SU(6)) \cdot \Z_2 \to (E_7)^{\gamma, \gamma'}$ satisfies
$$
    \sigma = \varphi_7((1, 1, F_{1, -1}), 1),  
    \sigma' = \varphi_7((1, 1, F_{-1, 1}), 1), \iota = \varphi_7((1, 1, F_{e_1}), 1)
$$
where} $F_{1, -1} = \diag(1, -1, -1, 1, -1, -1), F_{-1, 1} = \diag(-1, -1, 1, -1, -1, 1),F_{e_1}=$ $\diag(e_1, e_1,e_1,-e_1,-e_1,-e_1) \in SU(6).$
\vspace{2mm}

The group $\Z_2 = \{1, \gamma_1\}$ acts on the group $U(1)^{\times 7}$ by
$$
    \gamma_1(p, q, a_1, a_2, a_3, a_4, a_5) = (\ov{p}, \ov{q}, \ov{a}_4, \ov{a}_5, \ov{a}_1, \ov{a}_2, \ov{a}_3) $$
Let $(U(1)^{\times 7}) \cdot \Z_2$ be the semi-direct product of those groups under this action.
\vspace{2mm}

Now, we determine the structures of the group $(E_7)^{\gamma, \gamma', \sigma, \sigma',\iota} = ((E_7)^{\gamma, \gamma'})^{\sigma, \sigma',\iota}$. 
\vspace{2mm}

{\sc Theorem 4.3.} \qquad \qquad $((E_7)^{\gamma, \gamma', \sigma, \sigma',\iota})_0 
\cong U(1)^{\times 7}$. 
\vspace{2mm}

{\sc Proof.} For $\alpha \in (E_7)^{\gamma, \gamma', \sigma, \sigma',\iota} \subset 
(E_7)^{\gamma, \gamma'}$, there exist $p, q \in U(1)$ and $A \in SU(6)$ such 
that $\alpha = \varphi_7((p, q, A), 1)$ or $\alpha = \varphi_7((p, q, A), \gamma_1)$ 
(Proposition 4.1). For the case of $\alpha = \varphi_7((p, q, A), 1)$, by combining the conditions $\sigma\alpha\sigma 
= \alpha, \sigma'\alpha\sigma' = \alpha$ and $\iota\alpha\iota^{-1} = \alpha$ with  Lemma 4.2, we have
$$
    \varphi_7((p, q, F_{1, -1}AF_{1, -1}), 1) = \varphi_7((p, q, A), 1), \varphi_7((p, q, F_{-1, 1}AF_{-1, 1}), 1) = \varphi_7((p, q, A), 1)$$
\vspace{-1mm}
and
\vspace{-1mm}
$$
    \varphi_7((p, q, F_{e_1}A{F_{e_1}}^{-1}), 1) = \varphi_7((p, q, A), 1). $$
Hence
$$
\text{(i)}\,\,   F_{1, -1}AF_{1, -1} = A,
\,\,\,\,\text{(ii)}\left\{\begin{array}{l}
   p = \omega_1p
\vspace{1mm}\\
   q = \omega_1q
\vspace{1mm}\\
   F_{1, -1}AF_{1, -1} = \omega_1A,
\end{array} \right.
\,\,\,\,\text{(iii)}\left\{\begin{array}{l}
   p = {\omega_1}^2p
\vspace{1mm}\\
   q = {\omega_1}^2q
\vspace{1mm}\\
   F_{1, -1}AF_{1, -1} = {\omega_1}^2A,
\end{array} \right. 
$$
\vspace{-3mm}

$$
\text{(iv)}\,\,   F_{-1, 1}AF_{-1, 1} = A,
\,\,\,\,\text{(v)}\left\{\begin{array}{l}
   p = \omega_1p
\vspace{1mm}\\
   q = \omega_1q
\vspace{1mm}\\
   F_{-1, 1}AF_{-1, 1} = \omega_1A,
\end{array} \right.
\,\,\,\,\text{(vi)}\left\{\begin{array}{l}
   p = {\omega_1}^2p
\vspace{1mm}\\
   q = {\omega_1}^2q
\vspace{1mm}\\
   F_{-1, 1}AF_{-1, 1} = {\omega_1}^2A.
\end{array} \right. 
$$
and
$$
\text{(vii)}\,\,   F_{e_1}A{F_{e_1}}^{-1} = A,
\,\,\,\,\text{(viii)}\left\{\begin{array}{l}
   p = \omega_1p
\vspace{1mm}\\
   q = \omega_1q
\vspace{1mm}\\
   F_{e_1}A{F_{e_1}}^{-1} = \omega_1A,
\end{array} \right.
\,\,\,\,\text{(ix)}\left\{\begin{array}{l}
   p = {\omega_1}^2p
\vspace{1mm}\\
   q = {\omega_1}^2q
\vspace{1mm}\\
   F_{e_1}A{F_{e_1}}^{-1} = {\omega_1}^2A.
\end{array} \right. $$
We can eliminate the case (ii), (iii), (v), (vi), (viii) or (ix) because $p \ne 0$ or $q \ne 0$. Thus we have $p, q \in U(1)$ and $A \in S(U(1)^{\times 6})$. Since the mapping
$U(1)^{\times 5} \to S(U(1)^{\times 6})$, 
$$
h(a_1, a_2, a_3, a_4, a_5) = (a_1, a_2, a_3, a_4, a_5, \ov{a_1a_2a_3a_4a_5}) 
$$
is an isomorphism, the group satisfying with the conditions of case (i),(iv) and (vii) is $(U(1)^{\times 7})/\Z_3$. For the case of $\alpha=\varphi_7((p, q, A),\gamma_1)$, from $\varphi_7((p, q, A), \gamma_1)\!=\! \varphi_7((p, q, A),$ $ 1)\gamma_1, \varphi_7((1, 1, F_{1, -1}), 1)\gamma_1 = 
\gamma_1\varphi_7((1, 1, F_{1, -1}), 1), \varphi_7((1, 1, F_{-1, 1}), 1)\gamma_1= \gamma_1\varphi_7((1,$ $ 1, F_{-1, 1}), 1)$ and $\varphi_7((1, 1, F_{e_1}), 1)\gamma_1= \gamma_1\varphi_7((1,$ $ 1, F_{e_1}), 1)$, this case is in the same 
situation as above. Thus we have $(E_7)^{\gamma,\gamma',\sigma,\sigma', \iota}
\cong ((U(1)^{\times 7})/\Z_3) \cdot 
\Z_2,$ $ \Z_3 = \{(1)^{\times 7}, (w_1)^{\times 7}, ({w_1}^2)^{\times 7}\}$. The 
group $(U(1)^{\times 7})/\Z_3$ is naturally isomorphic to the torus 
$U(1)^{\times 7}$, 
hence we obtain $(E_7)^{\gamma, \gamma', \sigma, \sigma',\iota} 
\cong (U(1)^{\times 7}) \cdot \Z_2$. Therefore we have the required isomorphism of the theorem.   \hfill $\Box$
\vspace{3mm}

{\bf 3. The group $E_8$}
\vspace{2mm}

In the  $C$-vector space ${\gge_8}^C$:
$$
 {\gge_8}^C = {\gge_7}^C \oplus \gP^C \oplus \gP^C \oplus C \oplus C \oplus C,
$$
\noindent if we define the Lie bracket [$R_1, R_2$] by
$$
     [({\varPhi}_1, P_1, Q_1, r_1, u_1, v_1), ({\varPhi}_2, P_2, Q_2, r_2, u_2, v_2)] = ({\varPhi}, P, Q, r, u, v),
$$
$$
  \left \{ \begin{array}{l}
   {\varPhi}\, = [{\varPhi}_1, {\varPhi}_2] + P_1 \times Q_2 - P_2 \times Q_1 
\vspace{2mm}\\
   P = {\varPhi}_1P_2 - {\varPhi}_2P_1 + r_1P_2 - r_2P_1 + u_1Q_2 - u_2Q_1 
\vspace{2mm}\\
   Q = {\varPhi}_1Q_2 - {\varPhi}_2Q_1 - r_1Q_2 + r_2Q_1 + v_1P_2 - v_2P_1 
\vspace{2mm}\\
   r \,\,= - \displaystyle{\frac{1}{8}}\{P_1, Q_2\} + \displaystyle{\frac{1}{8}}\{P_2, Q_1\} + u_1v_2 - u_2v_1 
\vspace{2mm}\\
   u \,\,= \;\;\, \displaystyle{\frac{1}{4}}\{P_1, P_2\} + 2r_1u_2 - 2r_2u_1 
\vspace{2mm}\\
   v \,\,= - \displaystyle{\frac{1}{4}}\{Q_1, Q_2\} - 2r_1v_2 + 2r_2v_1, 
\end{array} \right. $$
\noindent then, ${\gge_8}^C$ becomes a simple $C$-Lie algebra of type ${E_8}$. \vspace{2mm}

The group ${E_8}^C$ is defined to be the automorphism group of the Lie algebra 
${\gge_8}^C$:
$$
{E_8}^C = \{\alpha \in \Iso_C({\gge_8}^C) \,| \, \alpha[R_1, R_2] = [\alpha R_1, \alpha R_2] \}. $$
                                                                                We define $C$-linear transformations $\sigma, {\sigma}', \wti{\lambda}$ of ${\gge_8}^C$ respectively by
$$
\begin{array}{c}
    \sigma({\varPhi}, P, Q, r, u, v) = (\sigma{\varPhi}\sigma, \sigma P, \sigma Q, r, u, v ),
\vspace{1mm}\\

    {\sigma}'({\varPhi}, P, Q, r, u, v) = ({\sigma}'{\varPhi}{\sigma}', {\sigma}' P, {\sigma}' Q, r, u, v ),
\vspace{1mm}\\
    \wti{\lambda}({\varPhi}, P, Q, r, u, v) = (\lambda{\varPhi}\lambda^{-1}, \lambda Q, -\lambda P,  -r, -v, -u),
\end{array} $$                                                                  \noindent where
$$ 
\begin{array}{c}
     \sigma{\varPhi}(\phi, A, B, \nu)\sigma 
               = {\varPhi}(\sigma\phi\sigma, \sigma A, \sigma B, \nu), 
\vspace{1mm}\\
      {\sigma}'{\varPhi}(\phi, A, B, \nu){\sigma}' =  { \varPhi}({\sigma}'\phi{\sigma}', {\sigma}' A, {\sigma}' B, \nu ),
\vspace{1mm}\\
         \lambda{\varPhi}(\phi, A, B, \nu)\lambda^{-1} = {\varPhi} (-{}^t\phi, -B, -A, -\nu).
\end{array} $$
\noindent ($\sigma, {\sigma}', \lambda$ of the left sides are the same ones used in [3].) Moreover, the complex conjugation in ${\gge_8}^C$ is denoted by $\tau$:
$$ 
       \tau({\varPhi}, P, Q, r, u, v) = (\tau{\varPhi}\tau, \tau P, \tau Q, 
                                          \tau r, \tau u, \tau v), $$
\noindent where $\tau{\varPhi}(\phi, A, B, \nu)\tau = {\varPhi}(\tau\phi\tau, \tau A, \tau B, \tau \nu).$
\vspace{2mm}

Now, we define the Lie group $E_8$ as a compact form of the complex  Lie group ${E_8}^C$ by
$$
        E_8 = \{ \alpha \in {E_8}^C\,|\,\tau\wti{\lambda}\alpha = \alpha\wti{\lambda}\tau \}.$$
Then, $E_8$ is a simply connected compact simple Lie group of type $E_8$. Note that $\sigma, {\sigma}', \wti{\lambda} \in E_8$. The Lie algebra $\gge_8$ of the Lie group $E_8$ is given by
\begin{eqnarray*}
        \gge_8 \!\!\! &=& \!\!\! \{R \in {\gge_8}^C \, |\, \tau\wti{\lambda}R = R \} 
\vspace{1mm}\\
         \!\!\! &=& \!\!\! \{({\varPhi}, P, -\tau\lambda P, r, u, -\tau u) \in {\gge_8}^C\,|\, {\varPhi} \in \gge_7, P \in \gP^C, r \in i\R, u \in C\}.
\end{eqnarray*}
Now, we will investigate the Lie algebra $(\gge_8)^{\sigma,\sigma'}$ of the group
$$
   (E_8)^{\sigma,\sigma'} = ((E_8)^{\sigma})^{\sigma'} = (E_8)^{\sigma} \cap (E_8)^{\sigma'}. $$ 

{\small \bibliographystyle{amsalpha}

\bigskip
\begin{flushright}

\begin{tabular}{l}
Toshikazu Miyashita \\
Komoro high school \\
Nagano, 384-0801, Japan \\
E-mail spin15ss16@ybb.ne.jp
\end{tabular}

\end{flushright}

\end{document}